\documentclass[psamsfonts, reqno]{amsart}

\usepackage{amscd,amsmath,amssymb,amsfonts,verbatim}
\usepackage{times}
\usepackage[cmtip, all]{xy}
\usepackage{graphicx}
\usepackage[active]{srcltx}
\usepackage{color}
\usepackage{tikz-cd}
\usepackage{textcomp}

\definecolor{refkey}{gray}{.5}   
\definecolor{labelkey}{gray}{.5} 
\definecolor{Red}{rgb}{1,0,0}

\newcommand{\pf}{{\bf Proof : }}

\newcommand{\qedwhite}{\hfill \ensuremath{\Box}}

\newtheorem{theo}{Theorem}[section]

\newtheorem{lem}[theo]{Lemma}

\theoremstyle{definition}

\newtheorem{defi}[theo]{Definition}

\title{Generalised Mennicke--Newman Lemma}
\author{Sampat Sharma}
\newcommand{\Addresses}{{
  \bigskip
  \footnotesize
  
  \textsc{Sampat Sharma, Department of Mathematics, IIT Bombay \\
            Mumbai 400076, INDIA}\par\nopagebreak
  \textit{E-mail:} Sampat Sharma \texttt{<sampat@math.iitb.ac.in; sampat.iiserm@gmail.com>}

  \medskip

  }}

\begin{document}
\maketitle
\subjclass 2010 Mathematics Subject Classification:{13C10, 13H99, 19B10, 19B14.}

 ~ {Mennicke--Newman lemma, Mennicke symbol, Unimodular rows.}
 \begin{abstract}
  Mennicke--Newman lemma for unimodular rows was used by W. van der Kallen to give a group structure on the orbit set $\frac{Um_{n}(R)}{E_{n}(R)}$ for a commutative noetherian ring of dimension $d\leq 2n-4.$ In this paper, we generalise the Mennicke--Newman lemma for $m\times n $ right invertible matrices.

 \end{abstract}
 
 \vskip0.50in

\begin{flushleft}
 Throughout this article we will assume $R$ to be a commutative ring with $1 \neq 0 .$
\end{flushleft}

\section{Introduction}
W. van der Kallen studied the orbit set $\frac{Um_{n}(R)}{E_{n}(R)}$ for a commutative noetherian ring $R.$ It was shown in \cite{vdk1} that the orbit set $\frac{Um_{n}(R)}{E_{n}(R)}$ carries a group structure for $n = d+1\geq 3,$ where $d$ is the 
dimension of the commutative noetherian ring $R.$ Mennicke--Newman lemma played an important role in defining the group structure on the orbit set $\frac{Um_{d+1}(R)}{E_{d+1}(R)}$. In \cite{vdk2}, W. van der Kallen extended the group structure to the 
range $2\leq d\leq 2n-4$ by establishing a bijection between the orbit set $\frac{Um_{n}(R)}{E_{n}(R)}$ and the weak Mennicke symbol $\mbox{WMS}_{n}(R),$ where $d$ denotes the stable dimension of $R.$
\par In this article we study the orbit set $\frac{Um_{m,n}(R)}{E_{n}(R)}$ for right invertible matrices $Um_{m,n}(R).$ We generalise the Mennicke--Newman lemma proved in \cite{vdk3} and prove the following : 
\begin{lem}
 Let $R$ be a commutative ring of $\mbox{sdim} ~ d$ and $S, T \in Um_{m,n}(R)$  such that $d\leq 2(n-m)-1$ , $ m\geq 2, n>m.$ Then there exists $\varepsilon_{1}, \varepsilon_{2} \in E_{n}(R)$ such that 
$$S\varepsilon_{1} = (X, \alpha), T\varepsilon_{2} = (I-X, \alpha)$$ 
where $X\in M_{m\times m}(R), \alpha \in M_{m\times (n-m)}(R).$
\end{lem}
\par The term $\frac{Um_{2,n}(R)}{E_{n}(R)}$ is just an orbit set of unimodular $2\times n$ matrices. We expect it to have a group structure for $n\geq 4, d\leq 2n-6.$ Indeed the analogy with algebraic topology in $(${\cite [Section 6]
{vdk3}}$)$ predicts an abelian group structure for $d\leq 2n-6.$ (Such analogy also gave correct predictions in \cite{vdk2} to get the group structure in certain range.) In this article we study the orbit set $\frac{Um_{2,n}(R[X])}{E_{n}(R[X])}$ for a local ring $R$ and 
prove that in some cases the map 
$$\mbox{Row}_{1}: \frac{Um_{2,n}(R[X])}{E_{n}(R[X])} \longrightarrow \frac{Um_{1,n}(R[X])}{E_{n}(R[X])}$$
is surjective.

\section{Preliminaries}
\par 
Let $v = (a_{0},a_{1},\ldots, a_{r}), w = (b_{0},b_{1},\ldots, b_{r})$ be two rows of length $r+1$ over a commutative ring 
$R$. A row $v\in R^{r+1}$ is said to be unimodular if there is a $w \in R^{r+1}$ with
$\langle v ,w\rangle = \Sigma_{i = 0}^{r} a_{i}b_{i} = 1$ and $Um_{r+1}(R)$ will denote
the set of unimodular rows (over $R$) of length 
$r+1$.
\par 
The group of elementary matrices,  denoted by $E_{r+1}(R)$,  is a subgroup of $GL_{r+1}(R)$ and is generated by the matrices 
of the form $e_{ij}(\lambda) = I_{r+1} + \lambda E_{ij}$, where $\lambda \in R, ~i\neq j, ~1\leq i,j\leq r+1,~
E_{ij} \in M_{r+1}(R)$ whose $ij^{th}$ entry is $1$ and all other entries are zero. The elementary linear group 
$E_{r+1}(R)$ acts on the rows of length $r+1$ by right multiplication. Moreover, this action takes unimodular rows to
unimodular
rows and we denote set of orbits of this action by $\frac{Um_{r+1}(R)}{E_{r+1}(R)}.$ 
The equivalence class of a row $v$ under this equivalence relation is denoted by $[v].$

\begin{defi}{ Mennicke symbol $\mbox{MS}_{n}(R), n\geq 2:$}
Following Suslin \cite{definedminnickesuslin}, we say that a Mennicke symbol of order $n$ on $R$ is a map $\phi$ 
from $Um_{n}(R)$ to an 
abelian group $G$ such that $\mbox{MS1}$ and $\mbox{MS2}$ holds:
\begin{align*}&\mbox{MS1}:~\mbox{ For~every~elementary~matrix~} \varepsilon \in
E_{n}(R) ~\mbox{and ~every}~v\in Um_{n}(R),\\~  
&~\mbox{we~have}~\phi(v\varepsilon) = \phi(v).\\
& \mbox{MS2}:~ \phi(x, a_{2},\ldots,a_{n})\phi(y,a_{2},\ldots,a_{n}) = \phi(xy,a_{2},\ldots,a_{n}).
\end{align*}
 
\end{defi}

\par In \cite{definedminnickesuslin}, it was proved by Suslin that there exists 
a universal Mennicke $n$-symbol denoted by $\mbox{ms}: Um_{n}(R) \longrightarrow \mbox{MS}_{n}(R).$ 

\begin{defi}{{Stable range condition}}~{{$\mbox{sr}_{n}(R):$}} We shall say stable range condition $\mbox{sr}_{n}(R)$ holds for $R$ if 
 for any $(a_{1},a_{2},\ldots,a_{n+1})\in Um_{n+1}(R)$ there exists $c_{i}\in R$ such that 
 $$(a_{1}+ c_{1}a_{n+1}, a_{2}+c_{2}a_{n+1},\ldots,a_{n}+c_{n}a_{n+1})\in Um_{n}(R).$$
 
\end{defi}
\begin{defi}{Stable range $\mbox{sr}(R),$ Stable dimension $\mbox{sdim}(R):$} We shall define the stable range
of $R$ denoted by $\mbox{sr}(R)$ to 
be the least integer $n$ such that $\mbox{sr}_{n}(R)$ holds. We shall define stable dimension of $R$ by 
$\mbox{sdim}(R) = \mbox{sr}(R) - 1.$
 
\end{defi}
\section{Mennicke--Newman Lemma for $Um_{n}(R)$}

\par
In Borsuk's cohomotopy groups one has the relation $\mbox{MS3}$:
\begin{align*}&\mbox{MS3}:~ \phi(x, a_{2},\ldots,a_{n})\phi(y,a_{2},\ldots,a_{n}) = \phi(xy,a_{2},\ldots,a_{n})~\mbox{if}~x+y=1.
 \end{align*}

\par
Now we list a few more relations that may or may not hold.
\begin{align*}
 &\mbox{MS4}:~ \phi(f^{2}, a_{2},\ldots,a_{n})\phi(g,a_{2},\ldots,a_{n}) = \phi(f^{2}g,a_{2},\ldots,a_{n}).\\
 &\mbox{MS5}:~ \phi(r, a_{2},\ldots,a_{n})\phi(1+q,a_{2},\ldots,a_{n}) = \phi(q,a_{2},\ldots,a_{n}),~\\
& \mbox{if} ~
 r(1+q) \equiv q~\mbox{mod}~(a_{2},\ldots,a_{n}).\\
 &\mbox{MS6}:~\phi(x, a_{2},\ldots,a_{n}) = \phi(-x,a_{2},\ldots,a_{n}).\\
 &\mbox{MS7}:~\phi(x, a_{2},\ldots,a_{n})^{m} = \phi(x^{m},a_{2},\ldots,a_{n})~\mbox{for}~m\geq 2.
\end{align*}

\begin{lem}
$(${\cite [Lemma 3.1]
{vdk3}}$)$
 Under $\mbox{MS1}$ and $\mbox{MS4},$ relations $\mbox{MS3}$ and $\mbox{MS5}$ are equivalent.
\end{lem}
${\pf}$ {\bf${\mbox{MS3} \Rightarrow \mbox{MS5}}:$} Since $r(1+q) \equiv q ~ \mbox{mod}(a_{2}, \ldots, a_{n}),$ we have 
$$r \equiv q(1-r)~\mbox{mod}(a_{2}, \ldots, a_{n}).$$
\begin{align*}
 \phi(q, a_{2}, \ldots, a_{n}) & = \phi(r(1+q), a_{2}, \ldots, a_{n})\\
 & = \phi(q(1-r)(1+q), a_{2}, \ldots, a_{n})\\
 & = \phi(r(1-r)(1+q)^{2}, a_{2}, \ldots, a_{n})\\
 & = \phi(r(1-r), a_{2}, \ldots, a_{n})\phi((1+q)^{2}, a_{2}, \ldots, a_{n})\\
 & = \phi(r, a_{2}, \ldots, a_{n})\phi((1-r), a_{2}, \ldots, a_{n})\phi((1+q)^{2}, a_{2}, \ldots, a_{n})\\
 & = \phi(r, a_{2}, \ldots, a_{n})\phi((1-r)(1+q)^{2}, a_{2}, \ldots, a_{n})\\
 & = \phi(r, a_{2}, \ldots, a_{n})\phi((1+q), a_{2}, \ldots, a_{n}).
\end{align*}
Last equality holds because $(1-r)(1+q) \equiv 1~\mbox{mod}(a_{2},\ldots, a_{n}).$
\par
{\bf${\mbox{MS5} \Rightarrow \mbox{MS3}}:$} Since $(y,a_{2},\ldots,a_{n})\in Um_{n}(R),$ there exists $v\in R$ such that $vy \equiv 1~\mbox{mod} (a_{2}, \ldots, a_{n}).$ Since 
$x + y = 1,$ we have 
$$vx \equiv (v-1)~ \mbox{mod} (a_{2}, \ldots, a_{n}).$$ 
\par 
In view of $\mbox{MS5}$ we have,
$$\phi(x, a_{2}, \ldots, a_{n})\phi(v, a_{2}, \ldots, a_{n}) = \phi(v-1, a_{2}, \ldots, a_{n}).$$
As $vy \equiv 1~\mbox{mod} (a_{2}, \ldots, a_{n}),$ we have 
$$\phi(x, a_{2}, \ldots, a_{n})\phi(v^{2}, a_{2}, \ldots, a_{n})\phi(y, a_{2}, \ldots, a_{n}) = \phi(v-1, a_{2}, \ldots, a_{n})
.$$
As $v-1 \equiv vx ~\mbox{mod} (a_{2}, \ldots, a_{n}),$ one has $ v-1 \equiv (v-1)vy \equiv v^{2}xy ~\mbox{mod} (a_{2}, \ldots, a_{n}).$ 
Therefore by $\mbox{MS1}$,  $$\phi(v-1, a_{2}, \ldots, a_{n}) =  \phi(v^{2}xy, a_{2}, \ldots, a_{n}).$$
Now, we have
$$\phi(x, a_{2}, \ldots, a_{n})\phi(y, a_{2}, \ldots, a_{n})\phi(v^{2}, a_{2}, \ldots, a_{n}) = 
\phi(xy, a_{2}, \ldots, a_{n})\phi(v^{2}, a_{2}, \ldots, a_{n}).$$
Thus 
$$\phi(x, a_{2},\ldots,a_{n})\phi(y,a_{2},\ldots,a_{n}) = \phi(xy,a_{2},\ldots,a_{n}).$$
 $~~~~~~~~~~~~~~~~~~~~~~~~~~~~~~~~~~~~~~~~~~~~~~~~~~~~~~~~~~~~~~~~~~~~~~~~~~~~~~~~~~~~~~~~~~~~~~~~~~~~~~~~~~~~~~\qedwhite$
 
\par We note Mennicke--Newman lemma proved by W. van der Kallen in $(${\cite[Lemma 3.2]{vdk3}}$),$  we write a proof here for 
completeness.
\begin{lem}
\label{mennickenewman}
 Let $R$ be a commutative ring with $\mbox{sdim}(R) \leq 2n-3.$ Let $v, w\in Um_{n}(R).$ Then there exists 
 $\varepsilon, \delta \in E_{n}(R)$ and $x, y, a_{i}\in R$ such that 
 $$v\varepsilon = (x, a_{2}, \ldots, a_{n}), w\delta = (y, a_{2}, \ldots, a_{n}), x+y = 1.$$
\end{lem}
${\pf}$ Let $v = (v_{1}, \ldots, v_{n}), w = (w_{1}, \ldots, w_{n}).$ We will prove the lemma in two steps, first we will show that using 
elementary transformations one can make sum of first cordinates of $v, w $ equals one. In next step using elementary transformations
we can make $v_{i} = w_{i}, 2\leq i\leq n.$
\par Since $  (v_{1}, \ldots, v_{n}), w = (w_{1}, \ldots, w_{n}) \in Um_{n}(R),$ we have $$(v_{1}w_{1}, v_{2}, \ldots, v_{n}, w_{2}, 
\ldots, w_{n})\in Um_{2n-1}(R).$$ Since $\mbox{sr}(R) = \mbox{sdim}(R) + 1 \leq 2n-2,$ we can add multiples of $v_{1}w_{1}$ to 
$v_{2}, \ldots, v_{n}, w_{2}, \ldots, w_{n}$ to make $(v_{2}, \ldots, v_{n}, w_{2}, \ldots, w_{n}) \in Um_{2n-2}(R),$ i.e. adding multiples of $v_{1}$ to $v_{2}, \ldots, v_{n}$ and multiples of $w_{1}$ to $w_{2}, \ldots, w_{n}$ we get $(v_{2}, \ldots, v_{n}, w_{2}, \ldots, w_{n}) \in Um_{2n-2}(R).$ Let $v_{2}', \ldots, v_{n}', w_{2}', \ldots, w_{n}' \in R$ such that 
$\sum_{i=2}^{n}v_{i}v_{i}' + \sum_{j=2}^{n}w_{j}w_{j}' = 1.$ Add $v_{i}'(1-v_{1}-w_{1})$-multiple of $i^{th}$ columns of $v$ to $v_{1}$ for $2\leq i\leq n$ and for $2\leq j\leq n$ add $w_{j}'(1-v_{1}-w_{1})$-multiple of $j^{th}$ columns 
 of $w$ to $w_{1}$ to get $v_{1} + w_{1} = 1.$
\par Add $(-w_{i}+v_{i})$-multiple of $w_{1}( = (1-v_{1}))$ to $w_{i}$ and $(-v_{i} + w_{i})$-multiple of $v_{1}$ to $v_{i},$ for $2\leq i\leq n,$ to get $v_{i} = w_{i}.$
$~~~~~~~~~~~~~~~~~~~~~~~~~~~~~~~~~~~~~~~~~~~~~~~~~~~~~~~~~~~~~~~~~~~~~~~~~~~~~~~~~~~~~~~~~~~~~~~~~~~~~~~~~~~~~~\qedwhite$

\section{Mennicke--Newman Lemma for $Um_{m,n}(R), n > m\geq 2$}
In this section we prove the generalised form of Mennicke--Newman lemma.
\begin{defi}
 An $\alpha \in M_{m \times n}(R)$ is said to be right invertible if  $\exists ~\beta \in M_{n \times m}(R)$ 
 such that $\alpha \beta = I_{m}.$ 
 We will denote set of all $ m\times n$ right invertible matrices by $Um_{m,n}(R).$
\end{defi}
We note a result a W. van der Kallen $(${\cite [Lemma 3.1]{vdk3}}$)$ which allows us to do the row operations in a right invertible matrix $M\in Um_{m,n}(R).$

\begin{lem} Let $R$ be a commutative ring, $m < n,$ $n\geq 3.$ The action by left multiplication of $E_{m}(R)$ on $\frac{Um_{m,n}(R)}{E_{n}(R)}$ is trivial.
\end{lem}
\begin{lem}
\label{genmennickenewman} Let $R$ be a commutative ring and $S, T \in Um_{m,n}(R)$  such that $\mbox{sdim}(R) =d\leq 2(n-m)-1$ , $ m\geq 2,  n>m.$ Then there exists $\varepsilon_{1}, \varepsilon_{2} \in E_{n}(R)$ such that 
$$S\varepsilon_{1} = (X, \alpha), T\varepsilon_{2} = (I-X, \alpha)$$ 
where $X\in M_{m\times m}(R), \alpha \in M_{m\times (n-m)}(R).$
\end{lem}
${\pf}$ We will prove it by induction on $m.$ Let $m = 2$ and $ S,T\in Um_{2,n}(R).$ In view of Lemma \ref{mennickenewman}, using elementary transformation we can make $s_{11} + t_{11} = 1$. Thus $S$ and $T$ looks like 
 $$S = \begin{bmatrix}
                 a & s_{12}&\ldots & s_{1n}\\
                  s_{21}& s_{22}& \ldots & s_{2n}\\
                \end{bmatrix},  T = \begin{bmatrix}
1-a & t_{12}&\ldots & t_{1n}\\
                  t_{21}& t_{22}& \ldots & t_{2n}\\
                 \end{bmatrix}.$$ 
Upon adding $(-s_{21}-t_{21})$-times first row of $S$ to the second row of $S$ and $(-s_{21}-t_{21})$-times first row of $T$ to the second row of $T$, we may assume that $S$ and $T$ looks like 
 $$S = \begin{bmatrix}
                 a & s_{12}&\ldots & s_{1n}\\
                  g& s_{22}& \ldots & s_{2n}\\
                \end{bmatrix},  T = \begin{bmatrix}
1-a & t_{12}&\ldots & t_{1n}\\
                  -g& t_{22}& \ldots & t_{2n}\\
                 \end{bmatrix}$$ for some $g\in R.$ (Note that other entries in second row of $S$ and $T$ will also change but by abuse of notations we are keeping them same.)
\par Since $ g, s_{22}, \ldots , s_{2n}$ and $-g, t_{22}, \ldots , t_{2n}$ are unimodular rows, we have $$(s_{23}, \ldots, s_{2n}, t_{23}, \ldots, t_{2n}, s_{22}t_{22}, g)\in Um_{2n-2}(R).$$
\par In view of $(${\cite [Corollary 9.4]
{7}}$)$, we may add multiples of $g$ to $s_{23}, \ldots, s_{2n}, t_{23}, \ldots, t_{2n}$ to make $\mbox{ht}(s_{23}, \ldots, s_{2n}, t_{23}, \ldots, t_{2n})\geq 2n-4\geq d+1.$ Thus  $(s_{23}, \ldots, s_{2n}, t_{23}, \ldots, t_{2n})\in Um_{2n-4}(R).$ 
Thus there exist  $s_{23}', \ldots, s_{2n}', t_{23}', \ldots, t_{2n}'\in R$ such that 
$$s_{23}s_{23}' + \cdots + s_{2n}s_{2n}'+ t_{23}t_{23}' + \cdots + t_{2n}t_{2n}' = 1.$$  
\par Multiply each $j^{th}$ ($3\leq j\leq n$) column of $S$ by $(1-t_{22}-s_{22})s_{2j}'$ and each $j^{th}$ ($3\leq j\leq n$) column of $T$ by $(1-t_{22}-s_{22})t_{2j}'$ and add to their 2nd column respectively. 
Thus {sum~of~$22^{th}$~entries~of~S~and~T} $ =  s_{22} + \sum_{j=3}^{n}(1-t_{22}-s_{22})s_{2j}'s_{2j} + t_{22} + \sum_{j=3}^{n}(1-t_{22}-s_{22})t_{2j}'t_{2j} = 1.$
\par We add $(-s_{12}-t_{12})$-multiple of first column to the second column in both $S$ and $T$ to get $s_{12} + t_{12} = 0.$ Observe that still one has $s_{22} + t_{22} = 1.$ Therefore $S$ and $T$ looks like 
$$S = \begin{bmatrix}
                 a & s_{12}& s_{13}&\ldots & s_{1n}\\
                  g& s_{22}&s_{23}& \ldots & s_{2n}\\
                \end{bmatrix},  T = \begin{bmatrix}
1-a & -s_{12}&t_{13}&\ldots & t_{1n}\\
                  -g&1- s_{22}&t_{23}& \ldots & t_{2n}\\
                 \end{bmatrix}.$$  
\par Next we will make later columns of $S$ and $T$ to be equal. Add $ \begin{bmatrix}
                 1-a \\
                  -g\\
                \end{bmatrix}(s_{13}-t_{13}) +  \begin{bmatrix}
                 -s_{12} \\
                 1- s_{22}\\
                \end{bmatrix}(s_{23}-t_{23})$ to the third column of $T$ and add  $ \begin{bmatrix}
                 a \\
                  g\\
                \end{bmatrix}(-s_{13}+t_{13}) +  \begin{bmatrix}
                 s_{12} \\
                  s_{22}\\
                \end{bmatrix}(-s_{23}+t_{23})$ to the third column of $S$ to get their third columns equal. Similarly other columns of $S$ and $T$ can be made equal. Thus $S$ and $T$ looks like 
$$S = \begin{bmatrix}
                 a & s_{12}& s_{13}&\ldots & s_{1n}\\
                  g& s_{22}&s_{23}& \ldots & s_{2n}\\
                \end{bmatrix},  T = \begin{bmatrix}
1-a & -s_{12}&s_{13}&\ldots & s_{1n}\\
                  -g&1- s_{22}&s_{23}& \ldots & s_{2n}\\
                 \end{bmatrix}.$$ 
Set $X = \begin{bmatrix}
                 a & s_{12}\\
                  g& s_{22}\\
                \end{bmatrix}$ to get the desired result in the case $m=2.$
 \par
 Now by induction hypothesis assume that the result holds for all $S, T\in Um_{k,n}(R), k< m.$ Let $S, T\in Um_{m,n}(R).$ By induction hypothesis we may assume that $S$ and $T$ looks like 
$$S = \begin{bmatrix}
                 Y &~& \beta \\
                  a_{m1}& a_{m2}&a_{m3}& \ldots & a_{mn}\\
                \end{bmatrix},  T = \begin{bmatrix}
I-Y &~& \beta \\
                  b_{m1}&b_{m2}&b_{m3}& \ldots & b_{mn}\\
                 \end{bmatrix},$$ 
for some $Y = (y_{ij})\in M_{(m-1,m-1)}(R)$ and $\beta \in M_{(m-1,n-(m-1))}(R).$ Multiply first row of $S$ by $(-b_{m1}-a_{m1})$ and add it to the $m^{th}$ row of $S$ and  multiply first row of $T$ by $(-b_{m1}-a_{m1})$ and add it to the $m^{th}$ row of $T$ 
to get $a_{m1} = -b_{m1}.$ After this we multiply $i^{th}$ row of $S$ by $(-b_{mi}-a_{mi})$ and add it to the $m^{th}$ row of $S$ and  multiply $i^{th}$ row of $T$ by $(-b_{mi}-a_{mi})$ and add it to the $m^{th}$ row of $T$ to make $a_{mi} = -b_{mi}, 2\leq i\leq m-1.$ 
(Note that during these operations $a_{mi}, b_{mi}$ will change but by the abuse of notation we are writing them as $a_{mi}, b_{mi}$ and they do satisfy $a_{mi} = -b_{mi}, 1\leq i\leq m-1.$) 
\par 
Since $ (a_{m1}, a_{m2},a_{m3}, \ldots , a_{mn})\in Um_{n}(R)$  and  the last row of $T$ \\
$  (-a_{m1},-a_{m2},-a_{m3}, \ldots, -a_{m,m-1}, b_{mm}, \ldots, b_{mn}) \in Um_{n}(R)$, we have\\ $(a_{mm}b_{mm}, a_{m,m+1}, \ldots, a_{mn}, b_{m,m+1}, \ldots, b_{mn}, a_{m1}, \ldots, a_{m,m-1})$ is also unimodular. 
\par Thus in view of $(${\cite [Corollary 9.4]
{7}}$)$, we may assume that $\mbox{ht}(a_{m,m+1}, \ldots, a_{mn}, b_{m,m+1}, \ldots, b_{mn})\geq 2(n-m)\geq d+1.$ Therefore we have\\ $(a_{m,m+1}, \ldots, a_{mn}, b_{m,m+1}, \ldots, b_{mn}) \in Um_{2(n-m)}(R).$ Thus there exist elements\\ $a_{m,m+1}', \ldots, a_{mn}', b_{m,m+1}', \ldots, b_{mn}'\in R$ such that 
$$ a_{m,m+1}a_{m,m+1}' + \cdots + a_{mn}a_{mn}' + b_{m,m+1}b_{m,m+1}'+ \cdots + b_{mn}b_{mn}' = 1.$$ 
\par Multiply $(m+i)^{th}, ~\mbox{for}~ i\geq 1,$ column of $S$ by $a_{m,m+i}'(1-a_{mm}-b_{mm})$ and add it to the $m^{th}$ column of $S$. Similarly multiply $(m+i)^{th},~\mbox{for}~ i\geq 1,$ column of $T$ by $b_{m,m+i}'(1-a_{mm}-b_{mm})$ and add it to the $m^{th}$ column of $T$ to get 
$a_{mm}+b_{mm} = 1.$
\par In the next step, add $(-s_{im}-t_{im})$-multiple of $i^{th}$ column to the $m^{th}$ column in both $S$ and $T$, for $1\leq i\leq m-1$ to get $s_{im} = -t_{im}.$ Here $s_{im}, t_{im}$ denotes the $im^{th}$ entry in the matrix $S$ and $T$ respectively for $1\leq i\leq m-1.$ Observe that still one has $a_{mm} + b_{mm} = 1.$ Thus $S$ and $T$ looks like 

$$S = \begin{bmatrix}
                 X &~& \alpha_{1} \\
                  \end{bmatrix},  T = \begin{bmatrix}
I-X &~& \alpha_{2} \\
                 \end{bmatrix}$$ 
for some $X\in M_{m,m}(R)$ and $\alpha_{1}, \alpha_{2}\in M_{m, n-m}(R).$ 
Let $$\alpha_{1} = \begin{bmatrix}
                 a_{1,m+1}&\ldots&a_{1n} \\
\vdots &~& \vdots\\
                  a_{m, m+1}&  \ldots & a_{mn}\\
                \end{bmatrix}, \alpha_{2} = \begin{bmatrix}
                 b_{1,m+1}&\ldots&b_{1n} \\
\vdots &~& \vdots\\
                  b_{m, m+1}&  \ldots & b_{mn}\\
                \end{bmatrix}.$$
By column operations we perform $$ \begin{bmatrix}
                 a_{1,m+i} \\
\vdots \\
                  a_{m, m+i}\\
                \end{bmatrix} + X\begin{bmatrix}
                 -a_{1,m+i}+b_{1,m+i} \\
\vdots \\
                  -a_{m, m+i}+b_{m,m+i}\\
                \end{bmatrix}~\mbox{in} ~S$$ 
and $$ \begin{bmatrix}
                 b_{1,m+i} \\
\vdots \\
                  b_{m, m+i}\\
                \end{bmatrix} +(I- X)\begin{bmatrix}
                 a_{1,m+i}-b_{1,m+i} \\
\vdots \\
                  a_{m, m+i}-b_{m,m+i}\\
                \end{bmatrix}~\mbox{in} ~T$$ 
for $1\leq i\leq n-m$, to make later columns in $S$ and $T$ equal. This completes the proof. 
$~~~~~~~~~~~~~~~~~~~~~~~~~~~~~~~~~~~~~~~~~~~~~~~~~~~~~~~~~~~~~~~~~~~~~~~~~~~~~~~~~~~~~~~~~~~~~~~~~~~~~~~~~~~~~~\qedwhite$

\section{Some results about the orbit $\frac{Um_{2,n}(R)}{E_{n}(R)}$} 
In \cite{vdk3}, W. van der Kallen has studied the orbit set  $\frac{Um_{2,n}(R)}{E_{n}(R)}$ . The map 
$$\mbox{Row}_{1}: \frac{Um_{2,n}(R)}{E_{n}(R)} \longrightarrow \frac{Um_{1,n}(R)}{E_{n}(R)}$$ associates to the orbit of a matrix the orbit of its first row. In this section, we will prove some results about the orbit set 
$\frac{Um_{2,n}(R[X])}{E_{n}(R[X])}$ for a local ring $R$. 
\begin{lem}
\label{d+1} Let $R$ be a local noetherian ring of dimension $d, d\geq 2, \frac{1}{d!}\in R.$ Let $n\geq d+1$, then the map 
$$\mbox{Row}_{1}: \frac{Um_{2,n}(R[X])}{E_{n}(R[X])} \longrightarrow \frac{Um_{1,n}(R[X])}{E_{n}(R[X])}$$
is surjective.
\end{lem}
${\pf}$ In view of $(${\cite [Remark 1.4.3]
{invent}}$)$, we may assume that $R$ is a reduced ring.  Let $n=d+1.$ In view of $(${\cite [Theorem 2.4]
{invent}}$)$, every unimodular row of length $d+1$ over $R[X]$ is completable to an invertible matrix of determinant 1, thus the map is surjective. Now let $n> d+1.$ Let $v = (v_{1}, \ldots, v_{d+1}, v_{d+2},\ldots, v_{n})\in Um_{n}(R[X]).$  In view of 
$(${\cite [Proposition 1.4.4]
{invent}}$)$, we can elementarily transform $v$ to a vector $u = (u_{1}(X), u_{2}(X), c_{3}, \ldots, c_{d+1}, c_{d+2}, \ldots, c_{n})\in Um_{n}(R[X])$ such that $c_{n}$ is a non-zero-divisor,
 $c_{i}\in R~\mbox{and}~ 3\leq i\leq n.$ 
\par  Let $I$ be an ideal of $R$ generated by $c_{d+2}, \ldots, c_{n}$.  Let ${-}$ denotes the ideal modulo $I[X].$ Since $R$ is a local ring of dimension  $d,$  $\mbox{dim}(\frac{R}{I})\leq d$ and it is a local ring. If $\mbox{dim}(\frac{R}{I}) = d,$ then in view of Rao's result $(${\cite [Theorem 2.4]
{invent}}$)$, we can complete $(\overset{-}{u_{1}(X)}, \overset{-}{u_{2}(X)},\overset{-}{ c_{3}}, \ldots, \overset{-}{c_{d+1}})$ to a matrix in $\overset{-}{R}[X]$ of determinant $1.$  If $\mbox{dim}(\frac{R}{I}) < d,$ then in view of $(${\cite[Theorem 7.2]{4}}$),$ we can complete $(\overset{-}{u_{1}(X)}, \overset{-}{u_{2}(X)},\overset{-}{ c_{3}}, \ldots, \overset{-}{c_{d+1}})$ to a matrix of determinant 1 (In this case we can actually complete it to an elementary matrix).
 Let us assume that matrix $\overset{-}{M}\in SL_{d+1}(\overset{-}{R}[X])$ be a completion of  $(\overset{-}{u_{1}(X)}, \overset{-}{u_{2}(X)},\overset{-}{ c_{3}}, \ldots, \overset{-}{c_{d+1}}).$ Lift the matrix $\overset{-}{M}$  to a matrix $M\in M_{d+1}(R[X])$ whose first row is given by the vector $  (u_{1}(X), u_{2}(X), c_{3}, \ldots, c_{d+1})$. As $\mbox{det}(\overset{-}{M}) = \overset{-}{1},$ we may choose $b_{d+2}, \ldots, b_{n} \in R[X]$ 
such that $\sum_{j=2}^{n-d-1}(c_{d+j}b_{d+j+1}-c_{d+j+1}b_{d+j}) = \mbox{det}(M)-1.$ Thus the following matrix is in $$  \begin{bmatrix}
                 u_{1}(X) &u_{2}(X)& c_{3}&\ldots &c_{d+1}&c_{d+2}&\ldots&c_{n} \\
m_{21}&m_{22}&m_{23}&\ldots&m_{2,d+1}&b_{d+2}&\ldots&b_{n}\\
                  \end{bmatrix}\in Um_{2,n}(R[X])$$ and the class of this matrix is the preimage of $[v] = [u].$ Thus the map is surjective. 
$~~~~~~~~~~~~~~~~~~~~~~~~~~~~~~~~~~~~~~~~~~~~~~~~~~~~~~~~~~~~~~~~~~~~~~~~~~~~~~~~~~~~~~~~~~~~~~~~~~~~~~~~~~~~~~\qedwhite$
 
\begin{lem}
\label{dimension3}  Let $R$ be a local noetherian ring of dimension $3,  \frac{1}{2}\in R.$ Let $n\geq 3$, then the map 
$$\mbox{Row}_{1}: \frac{Um_{2,n}(R[X])}{E_{n}(R[X])} \longrightarrow \frac{Um_{1,n}(R[X])}{E_{n}(R[X])}$$
is surjective.
\end{lem}${\pf}$ In view of $(${\cite [Remark 1.4.3]
{invent}}$)$, we may assume that $R$ is a reduced ring. Let $n=3.$ Then in view of $(${\cite [Theorem 3.1]
{trans}}$)$ every unimodular row of length $3$ over $R[X]$ is completable to an invertible matrix of determinant 1, thus the map is surjective. Now let $n> 3.$ Let $v = (v_{1}, \ldots, v_{3},\ldots, v_{n})\in Um_{n}(R[X]).$ In view of 
$(${\cite [Proposition 1.4.4]
{invent}}$)$, we can elementarily transform $v$ to a vector $u = (u_{1}(X), u_{2}(X), c_{3}, \ldots,  \ldots, c_{n})\in Um_{n}(R[X])$ such that $c_{i}\in R, 3\leq i\leq n$ and $c_{3}$ is a non-zero-divisor.
\par  Let $I$ be an ideal of $R$ generated by $c_{4}, \ldots, c_{n}$.  Let ${-}$ denotes the ideal modulo $I[X].$ Since $\mbox{dim}(\frac{R}{I})\leq 3,$ we can complete $(\overset{-}{u_{1}(X)}, \overset{-}{u_{2}(X)},\overset{-}{ c_{3}}, )$ 
to a matrix $\overset{-}{M}$. Lift the matrix $\overset{-}{M}$  to a matrix $M$ whose first row is given by the vector $  (u_{1}(X), u_{2}(X), c_{3})$. As $\mbox{det}(\overset{-}{M}) = \overset{-}{1},$ we may choose $b_{4}, \ldots, b_{n} \in R[X]$ 
such that $\sum_{j=2}^{(n-1)/2}(c_{2i}b_{2i+1}-c_{2i+1}b_{2i}) = \mbox{det}(M)-1.$ Thus the following matrix is in $$Um_{2,n}(R[X]) :  \begin{bmatrix}
                 u_{1}(X) &u_{2}(X)& c_{3}&c_{4} &\ldots&c_{n} \\
m_{21}&m_{22}&m_{23}&b_{4}&\ldots&b_{n}\\
                  \end{bmatrix}$$ and the class of this matrix is the preimage of $[v] = [u].$ Thus the map is surjective. 
$~~~~~~~~~~~~~~~~~~~~~~~~~~~~~~~~~~~~~~~~~~~~~~~~~~~~~~~~~~~~~~~~~~~~~~~~~~~~~~~~~~~~~~~~~~~~~~~~~~~~~~~~~~~~~~\qedwhite$

\section{Some known results}
We note some known results in this section:

\begin{lem}$(${\cite [Theorem 7.1] {7}}$)$

Let $R$ be a commutative ring with unity with finite stable rank. Then $sr(I) \leq sr(R)$ and $sr(R/I) \leq sr(R)$ for any ideal $I$ in $R$.

\end{lem}

\begin{lem}$(${\cite [Lemma 7.5] {7}}$)$

Let $R$ be a commutative ring with unity. Suppose that the maximal ideal space of $R$ is Noetherian and a finite union of spaces of dimensions that do not exceed $m.$ Then the stable range of $R$ is less than or equal to $m + 1.$

\end{lem}

\begin{lem}$(${\cite [Lemma 7.6] {7}}$)$

Let $R$ be a commutative ring with unity, $m$ be a natural number, $I_1, I_2, \ldots, I_m$ be ideals such that $I_{i} \subset rad(R).$ Then, $sr(R) = \mbox{max}_{1\leq i\leq m} sr(R/I_{i}).$

\end{lem}

\begin{theo}$(${\cite [Theorem 16.4] {7}}$)$ Let $R$ be a finitely generated ring of dimension $1.$ Assume that no quotient ring of $R$ is totally imaginary arithmetic ring. Then, $SK_{1}(R) = 0.$

\end{theo}

\begin{theo}$(${\cite [Theorem 17.2] {7}}$)$ Let $R$ be a finitely generated ring. Assume that no quotient ring of $R$ is totally imaginary arithmetic ring (i.e. $kR = R$ for some natural number $k\geq 2.$) Then $sr(R) \leq \mbox{max}(2, dim(R)).$

\end{theo}

\begin{theo}$(${\cite [Corollary 17.3] {7}}$)$ Let $F \rightarrow R$ be a finitely generated algebra over a field $F,$ where $F$ is algebraic over a finite field. Then, $sr(R) \leq \mbox{max}(2, dim(R)).$

\end{theo}

\begin{theo}$(${\cite [Theorem 18.2] {7}}$)$ For any finitely generated ring $R,$ we have $sr(R) \leq \mbox{max}(3, dim(R)).$

\end{theo}

\begin{theo}$(${\cite [Corollary 19.1] {7}}$)$  Let $R$ be a ring without zero divisors, in which $0$ is not representable as sums of non-zero squares. Then $sr(R[X_1, X_2, \ldots,X_n]) \geq n+1$ for all natural numbers $n.$

\end{theo}

\begin{theo}$(${\cite [Theorem 19.3] {7}}$)$ For any $R$ of characteristic $0$ and natural number $n$, the stable range $sr(R[X_1, X_2, \ldots, X_n]) \geq (n+1)/2.$

\end{theo}

\begin{lem} $(${\cite [Corollary 3.5] {improved}}$)$ Let A be a regular affine algebra of Krull dimension 3 over
a C1 field k which is perfect if its characteristic is 2 or 3. Then the Vaser-
stein symbol $V : Um_{3} (A)/E_{3} (A) \rightarrow W_{E} (A)$ is an isomorphism. 
\end{lem}

\begin{lem}$(${\cite [Proposition 4.2] {improved}}$)$ The ring $P_3$ is a 3 dimensional ring for 'which the Vaserstein
symbol
$$V: \frac{Um_{3}(P_{3})}{E_{3}(P_{3})} \to W_{E}(P_{3})$$
is not injective. 
\end{lem}

\Addresses
\end{document}